\documentclass[12pt,leqno]{amsart}
\usepackage{amsmath, amssymb, latexsym, amsthm}
\usepackage[all]{xy}

\newcommand{\arx}[1]{\texttt{http://arxiv.org/abs/#1}}

\newcommand{\gr}{\mbox{\textit{\tiny gp}}}

\newcommand{\inv}{^{-1}}
\newcommand{\Cantor}{{\{0,1\}^\N}}

\newcommand{\seq}[1]{\{#1\}_{n\in\N}}

\newcommand{\cF}{\mathcal{F}}

\newcommand{\N}{\mathbb{N}}
\newcommand{\Null}{\mathcal{N}}
\newcommand{\NN}{{{}^{\N}\N}}
\newcommand{\NNup}{{{}^{{\N}\nearrow}\N}}
\renewcommand{\inf}{P_\oo(\N)}

\newcommand{\cU}{\mathcal{U}}
\newcommand{\cP}{\mathcal{P}}
\newcommand{\Union}{\bigcup}
\newcommand{\cV}{\mathcal{V}}
\newcommand{\cW}{\mathcal{W}}

\newcommand{\Impl}{\Rightarrow}
\long\def\forget#1\forgotten{}

\renewcommand{\b}{\mathfrak{b}}

\renewcommand{\i}{\item}

\newcommand{\oo}{\infty}

\newcommand{\Iff}{\Leftrightarrow}

\newcommand{\nin}{\not\in}

\newcommand{\sbst}{\subseteq}

\newcommand{\sm}{\setminus}

\newcommand{\add}{\mathsf{add}}

\newtheorem{thm}{Theorem} 

\newtheorem{lem}[thm]{Lemma}
\newtheorem{cor}[thm]{Corollary}

\theoremstyle{definition}

\theoremstyle{remark}
\newtheorem{rem}[thm]{Remark}
\newcommand{\be}{\begin{enumerate}}
\newcommand{\ee}{\end{enumerate}}
\newcommand{\bi}{\begin{itemize}}
\newcommand{\ei}{\end{itemize}}

\title[The Hurewicz property and slaloms]{The Hurewicz covering property and slaloms in the Baire space}
\author{Boaz Tsaban}
\thanks{Partially supported by the Golda Meir Fund and
the Edmund Landau Center for Research in Mathematical Analysis and Related Areas,
sponsored by the Minerva Foundation (Germany).}
\address{Einstein Institute of Mathematics, Hebrew University of Jerusalem,
Givat Ram, Jerusalem 91904, Israel}
\email{tsaban@math.huji.ac.il}
\urladdr{http://www.cs.biu.ac.il/\~{}tsaban}

\begin{document}
\begin{abstract}
According to a result of Ko\v{c}inac and Scheepers, the Hurewicz
covering property is equivalent to a somewhat simpler selection
property: For each sequence of large open covers of the space one
can choose finitely many elements from each cover to obtain a
groupable cover of the space. We simplify the characterization
further by omitting the need to consider sequences of covers: A
set of reals $X$ has the Hurewicz property if, and only if,
each large open cover of $X$ contains a groupable subcover.
This solves in the affirmative a problem of Scheepers.

The proof uses a rigorously justified abuse of notation and a
``structure'' counterpart of a combinatorial
characterization, in terms of slaloms,
of the minimal cardinality $\b$ of an unbounded family
of functions in the Baire space.
In particular, we obtain a new characterization of $\b$.
\end{abstract}

\maketitle

\section{Introduction}

A separable zero-dimensional metrizable space $X$ has the
\emph{Hurewicz property} \cite{HURE27} if:
\begin{quote}
For each sequence $\seq{\cU_n}$ of open covers of $X$
there exist finite subsets $\cF_n\sbst\cU_n$, $n\in\N$, such
that $X\sbst\bigcup_n\bigcap_{m>n} \cup\cF_m$.
\end{quote}
This property is a generalization of $\sigma$-compactness.

Much effort was put in the past in order to find a simpler characterization of this property.
In particular, it was desired to avoid the need to glue the elements of
each $\cF_n$ together (that is, by taking their union) in the definition
of the Hurewicz property.

The first step toward simplification was
the observation that one may restrict attention to sequences of \emph{large} (rather than
arbitrary) open covers of $X$ in the above definition \cite{coc1}
($\cU$ is a \emph{large cover} of $X$ if each member of $X$ is contained in infinitely many
members of $\cU$).

The main ingredient in the next major step toward this goal was
implicitly studied in \cite{smzpow, NSW, FunRez} while considering
spaces having analogous properties in all finite powers and a
close relative of the Reznichenko (or: weak Fr\'echet-Urysohn)
property. Finally, this ingredient was isolated and analyzed in
\cite{coc7}: A large cover $\cU$ of $X$ is \emph{groupable} if
there exists a partition $\cP$ of $\cU$ into finite sets such that
for each $x\in X$ and all but finitely many $\cF\in\cP$,
$x\in\cup\cF$. Observe that ignoring all but countably many
elements of the partition, we see that each groupable cover
contains a countable groupable cover.
Moreover, in \cite{splittability} it is proved that for the types
of spaces considered here, each large open cover contains a countable
large cover.

One of the main results in Ko\v{c}inac-Scheepers' \cite{coc7} is that
the Hurewicz property is equivalent to the following one:
\begin{quote}
For each sequence $\seq{\cU_n}$ of large open covers of $X$
there exist finite subsets $\cF_n\sbst\cU_n$, $n\in\N$, such
that $\bigcup_n \cF_n$ is a groupable cover of $X$.
\end{quote}
This characterization is misleading in its pretending to be a mere
result of unstitching the unions $\cup\cF_n$: The sets $\cF_n$
need not be disjoint in the original definition, and overcoming this
difficulty requires a deep analysis involving infinite game-theoretic
methods -- see \cite{coc7}.

In this paper we take the task of simplification one
step further by removing the need to consider \emph{sequences} of covers.
We prove that the Hurewicz property is equivalent to:
\begin{quote}
$(\star)$ \quad Each large cover of $X$ contains a groupable subcover.
\end{quote}
This solves in the affirmative a problem of Scheepers \cite[Problem 1]{smzpow},
which asks whether, for strong measure zero sets,
$(\star)$ is equivalent to the Hurewicz property.

Another way to view this simplification is as follows:
The Ko\v{c}inac-Scheepers characterization is
equivalent to requiring that the resulting cover $\bigcup_n \cF_n$ is large
\emph{together} with the property that each large open cover of $X$ contains a groupable cover of $X$.
The first requirement has also appeared in the literature, and is equivalent to a
property introduced by Menger in \cite{MENGER} (see \cite{coc1}).
Our result says that it is enough to require only that the second property holds,
or in other words, that the second property actually implies the first.

\section{Two possible interpretations}
Our proof relies on a delicate interplay between two possible interpretations
of the term ``large open cover of $X$'' when $X$ is a subspace of another
space $Y$:
\be
\i A large open cover of $X$ by subsets of $X$ which are relatively open in $X$; and
\i A large open cover of $X$ by open subsets of $Y$.
\ee
The notions do not coincide, because a large cover of the second type, when restricted
to $X$, need not be large -- it can even be finite.

For brevity, we will use the following notation.
For a space $X$, denote the property that
each large cover of $X$ by open subsets of $X$ contains a groupable cover
of $X$ by $\binom{\Lambda_X}{\Lambda^{\gr}_X}$.
We write $\binom{\Lambda}{\Lambda^{\gr}}$ instead of
$\binom{\Lambda_X}{\Lambda^{\gr}_X}$ when the space $X$ is
clear from the context. It is easy to see that
$\binom{\Lambda_X}{\Lambda^{\gr}_X}$ implies that each \emph{countable}
large open cover of $X$ is groupable (divide the countably many
remaining elements between the sets in the partition so that
they remain finite).

We will need the following simple fact.

\begin{lem}\label{preservation}
The property $\binom{\Lambda}{\Lambda^{\gr}}$ is preserved
under taking closed subsets and continuous images, that is:
\be
\i If $\binom{\Lambda_X}{\Lambda^{\gr}_X}$ holds
and $C$ is a closed subset of $X$, then
$\binom{\Lambda_C}{\Lambda^{\gr}_C}$ holds.
\i If $\binom{\Lambda_X}{\Lambda^{\gr}_X}$ holds and
$Y$ is a continuous image of $X$,
then $\binom{\Lambda_Y}{\Lambda^{\gr}_Y}$ holds.
\ee
\end{lem}
\begin{proof}
(1) Assume that $\cU$ is a large open cover of $C$.
Then $\tilde\cU=\{U\cup (X\sm C) : U\in\cU\}$ is
a large open cover of $X$. Applying the groupability of
$\tilde\cU$ for $X$ and forgetting the $X\sm C$ part of the open sets
shows the groupability of $\cU$ for $C$.

(2) Assume that $f:X\to Y$ is a continuous surjection and that
$\cU$ is a large open cover of $Y$ (by open subsets of $Y$).
Then $\cV=\{f\inv[U]:U\in\cU\}$ is a large open cover of $X$.
By the assumption, there exists a groupable subcover $\cW\sbst\cV$ for $X$.
It follows that $\{U\in\cU : f\inv[U]\in\cW\}$ is a groupable cover of $Y$.
\end{proof}

The following theorem tells us that \emph{for our purposes},
it does not matter which notion of large covers we use
(so that we can switch between the two notions at our convenience).
\begin{thm}\label{outercovers}
Assume that $X$ is a subspace of $Y$ and $\binom{\Lambda_X}{\Lambda^{\gr}_X}$ holds.
Then each countable collection $\cU$ of open sets in $Y$ which is a large cover of
$X$ is groupable for $X$.
\end{thm}
\begin{proof}
We will repeatedly use the following lemma.
\begin{lem}\label{onestep}
Assume that $X$ is a subspace of $Y$,
and $\cU=\seq{U_n}$ is a large open cover of $X$ by open subsets of $Y$.
Define an equivalence relation $\sim$ on $\N$ by
$$n\sim m\quad\mbox{if}\qquad X\cap U_n = X\cap U_m.$$
Let $A = \{n : [n]\mbox{ is infinite}\}$, and $V=\Union_{n\in A} U_n$.
Then $\{U_n : n\in A\}$ is a groupable cover of $X\cap V$,
and $\{U_n : n\nin A\}$ is a large cover of $X\sm V$
(by open subsets of $Y$).
\end{lem}
\begin{proof}
Define a partition of $A$ as follows:
Let $[n_0], [n_1], \dots$ enumerate the elements of $A/\sim$.
Let $F_0$ contain the first element of $[n_0]$,
$F_1$ contain the second element of $[n_0]$ and the first element of $[n_1]$,
$F_2$ contain the third element of $[n_0]$, the second element of $[n_1]$ and the first element of $[n_2]$,
etc.
Fix $x\in X\cap V$, and let $i$ be such that $x\in U_{n_i}$.
Then for all but finitely many $n$, there exists $k\in F_n\cap [n_i]$ and therefore $x\in U_k$.
This proves the first assertion.

Assume that $X\not\sbst V$.
As $\cU$ is a large cover of $X\sm V$ and for $x\in X\sm V$ and $n\in A$, $x\nin U_n$,
there must exist infinitely many $n\nin A$ such that $x\in U_n$.
\end{proof}
We now prove Theorem~\ref{outercovers}.
Enumerate $\cU$ bijectively as $\seq{U_n}$.
We make the following definition by transfinite induction on $\alpha<\aleph_1$
(and make sure that indeed it terminates at some $\alpha<\aleph_1$).
Carry out the following construction as long as $A_\alpha$ is not empty.
\be
\i \emph{First step:}
Set $X_0=X$, $B_0 = \N$, and $V_0 = \emptyset$.

\i \emph{Successor step:} Assume that $X_\alpha$, $B_\alpha$, and $V_\alpha$
are defined, and $\{U_n : n\in B_\alpha\}$ is a large cover of $X_\alpha\sm V_\alpha$.
Set $X_{\alpha+1}=X_\alpha\sm V_\alpha$,
and define an equivalence relation $\sim_{\alpha+1}$ on $B_\alpha$ by
$n\sim_{\alpha+1} m$ if $X_{\alpha+1}\cap U_n = X_{\alpha+1}\cap U_m$.
Let $A_{\alpha+1} = \{n\in B_\alpha : [n]_{\sim_{\alpha+1}}\mbox{ is infinite}\}$,
$B_{\alpha+1} = B_\alpha\sm A_{\alpha+1}$,
and $V_{\alpha+1} = \Union_{n\in A_{\alpha+1}} U_n$.
Use Lemma~\ref{onestep} to obtain a partition $\seq{F^{\alpha+1}_n}$ of $A_{\alpha+1}$ into finite sets
witnessing that $\{U_n : n\in A_{\alpha+1}\}$ is a groupable cover of $X_{\alpha+1}\cap V_{\alpha+1}$
(and $\{U_n : n\in B_{\alpha+1}\}$ is a large cover of $X_{\alpha+1}\sm V_{\alpha+1}$).

\i \emph{Limit step:}
Assume that $\alpha$ is a limit and the construction was carried up to step $\alpha$.
Set $X_\alpha = \bigcap_{\beta<\alpha}X_\beta$,
$A_\alpha = \bigcup_{\beta<\alpha} A_\beta$,
$B_\alpha = \bigcap_{\beta<\alpha} B_\beta = \N\sm A_\alpha$,
and
$V_\alpha = \bigcup_{\beta<\alpha} V_\beta$.
For each $x\in X_\alpha$ and each $\beta<\alpha$, $x\nin V_\beta$,
that is, $\{n : x\in U_n\}$ is disjoint from $A_\beta$.
Thus, $\{n : x\in U_n\}$ is infinite, and is a subset of $B_\alpha$.
In other words, $\{U_n : n\in B_\alpha\}$ is a large cover of $X_\alpha$.
Observe that in this case, $X_\alpha$ is disjoint from $V_\alpha$.
\ee
As long as the construction continues, $A_\alpha$ is not empty and therefore
$B_{\alpha+1}$ is a proper subset of $B_\alpha$. Thus, as $B_0\sbst\N$, the construction
cannot continue uncountably many steps. Let $\alpha<\aleph_1$ be the step where the
construction terminates (this can only happen when $\alpha$ is a successor).
Then $A_\alpha$ is empty, therefore $V_\alpha$ is empty,
thus $\{U_n : n\in B_\alpha\}$ is a large cover of $X_\alpha$.
The definition of $B_\alpha$ implies that in this case,
$\{U_n\cap X_\alpha : n\in B_\alpha\}$ is a large cover of $X_\alpha$ by open subsets of $X_\alpha$.
By the construction, $X_\alpha$ is a closed subset of $X$ (an intersection of closed sets).
By Lemma~\ref{preservation}, $\{U_n\cap X_\alpha : n\in B_\alpha\}$ is groupable for
$X_\alpha$; let $\seq{F^{\alpha+1}_n}$ be a partition of $B_\alpha$ into finite sets that
witnesses that.

The partitions $\seq{F_n^{\beta+1}}$ where $\beta\le\alpha$ form a countable
family of partitions of disjoint subsets of $\N$.
Relabel these partitions as $\{\seq{G^m_n} : m\in\N\}$, and define a partition
$\seq{H_n}$ of $\N$ into finite sets by
$$H_n = \Union_{\max\{i,j\} = n}G^i_j.$$
Observe that $X\sbst X_\alpha\cup\Union_{\beta<\alpha}(V_{\beta+1}\sm V_\beta)$,
where each $X\cap(V_{\beta+1}\sm V_\beta)$ is taken care by $\seq{F_n^{\beta+1}}$,
and $X_\alpha$ is taken care by $\seq{F_n^{\alpha+1}}$.
Consequently, for each $x\in X$ there exists $m$ such that
$x\in\cup\{U_k : k\in G^m_n\}\sbst\cup\{U_k : k\in H_n\}$
for all but finitely many $n$. This shows that $\cU$ is a groupable cover of $X$.
\end{proof}

\section{The main theorem}

\begin{thm}\label{HureChar}
For a separable and zero-dimensional metrizable space $X$, the following are equivalent:
\be
\i $X$ has the Hurewicz property,
\i Every large open cover of $X$ contains a groupable cover of $X$; and
\i Every countable large open cover of $X$ is groupable.
\ee
\end{thm}
\begin{proof}
$(2\Iff 3)$ By Proposition 1.1 of \cite{splittability}, every large open cover
of $X$ contains a \emph{countable} large open cover of $X$.

$(1\Impl 3)$ This is proved in \cite[Lemma 3]{smzpow} and \cite[Lemma 8]{coc7}.

$(3\Impl 1)$ We will prove this assertion by a sequence of small steps,
using the results of the previous section.

The \emph{Baire space} $\NN$ of infinite sequences
of natural numbers is equipped with the product topology
(where the topology of $\N$ is discrete).
A quasiordering $\le^*$ is defined on the Baire space $\NN$ by eventual
dominance:
$$f\le^* g\quad\mbox{if}\qquad f(n)\le g(n)\mbox{ for all but finitely many }n.$$
We say that a subset $Y$ of $\NN$ is \emph{bounded} if
there exists $g$ in $\NN$ such that for each $f\in Y$, $f\le^* g$.
Otherwise, we say that $Y$ is \emph{unbounded}.
According to a theorem of Hurewicz
\cite{HURE27} (see also Rec\l{}aw \cite{RECLAW}),
$X$ has the Hurewicz property if, and only if,
each continuous image of $X$ in $\NN$ is bounded.
Let $\NNup$ denote the subspace of $\NN$ consisting
of the strictly increasing elements of $\NN$.
The mapping from $\NN$ to $\NNup$ defined by
$$f(n)\mapsto g(n)=f(0)+\dots+f(n)+n$$
is a homeomorphism which preserves boundedness
in both directions.
Consequently, Hurewicz' theorem can be stated
using $\NNup$ instead of $\NN$.

For $f,g\in\NNup$, we say that
\emph{$f$ goes through the slalom defined by $g$}
if for all but finitely many $n$, there exists $m$ such that $f(m)\in [g(n),g(n+1))$.
A subset $Y$ of $\NNup$ \emph{admits a slalom} if there exists $g\in\NNup$
such that each $f\in Y$ goes through the slalom $g$.
\begin{lem}[folklore]\label{bddslalom}
Assume that $Y\sbst\NNup$. The following are equivalent:
\be
\i $Y$ is bounded,
\i $Y$ admits a slalom; and
\i There exists a partition $\seq{F_n}$ of $\N$ into finite sets
such that for each $f\in Y$
and all but finitely many $n$, there exists $m$ such that $f(m)\in F_n$.
\ee
\end{lem}
For completeness, we give a short proof.
\begin{proof}
$(1\Impl 2)$ Assume that $g\in\NNup$ bounds $Y$.
Define inductively $h\in\NNup$ by
\begin{eqnarray*}
h(0) & = & g(0)\\
h(n+1) & = & g(h(n))+1
\end{eqnarray*}
Then for each $f\in Y$ and all but finitely many $n$,
$h(n)\le f(h(n))\le g(h(n))<h(n+1)$, that
is, $f(h(n))\in [h(n),h(n+1))$.

$(2\Impl 1)$ Assume that $Y$ admits a slalom $g$.
Let $h$ be a function which eventually dominates all functions of the
form $f(n)=g(n_0+n)$, $n_0\in\N$.
Let $f$ be any element of $Y$ and choose $n_0$ such that for each
$n\ge n_0$, there exists $m$ such that $f(m)\in[g(n),g(n+1))$.
Choose $m_0$ such that $f(m_0)\in[g(n_0),g(n_0+1))$.
By induction on $n$, we have that $(f(n)\le) f(m_0+n)\le g(n_0+1+n)$
for all $n$. For large enough $n$, we have that $g(n_0+1+n)\le h(n)$,
thus $f\le^* h$.

Clearly $(2\Impl 3)$. We will show that $(3\Impl 2)$.
Let $\seq{F_n}$ be as in (2).
Define $g\in\NNup$ as follows:
Set $g(0)=0$.
Having defined $g(0),\dots,g(n-1)$,
let $m$ be the minimal such that $F_m\cap [0,g(n-1))=\emptyset$,
and set $g(n) = \max F_m +1$.
Then for each $n$ there exists $F_m$ such that $F_m\sbst[g(n),g(n+1))$.
Consequently, $Y$ admits the slalom defined by $g$.
\end{proof}

The \emph{Cantor space} $\Cantor$ is also equipped with the product topology.
Identify $P(\N)$ with $\Cantor$ by characteristic functions.
The \emph{Rothberger space} $\inf$ is the subspace of $P(\N)$ consisting of all
infinite sets of natural numbers.
The space $\NNup$ is homeomorphic to $\inf$ by
identifying each $f\in\NNup$ with its image $f[\N]$
(so that $f$ is the increasing enumeration of $f[\N]$).

Translating Lemma~\ref{bddslalom} into the language of
$\inf$ and using Hurewicz' theorem, we obtain the
following characterization of
the Hurewicz property in terms of continuous images
in the Rothberger space.

\begin{lem}\label{combgp}
For a separable and zero-dimensional metrizable space $X$, the following are equivalent:
\be
\i $X$ has the Hurewicz property,
\i For each continuous image $Y$ of $X$ in $\inf$
there exists $g\in\NNup$ such that
for each $y\in Y$, $y\cap [g(n),g(n+1))\neq\emptyset$
for all but finitely many $n$; and
\i For each continuous image $Y$ of $X$ in $\inf$
there exists a partition $\seq{F_n}$ of $\N$ into
finite sets such that for each $y\in Y$,
$y\cap F_n\neq\emptyset$ for all but finitely many $n$.
\ee
\end{lem}

Assume that every countable large open cover of $X$ is groupable.
We will show that (3) of Lemma~\ref{combgp} holds.
Let $Y$ be a continuous image of $X$ in $\inf$.
Then by Lemma~\ref{preservation}, $\binom{\Lambda_Y}{\Lambda^{\gr}_Y}$ holds.
Thus, by Theorem~\ref{outercovers},
every countable large open cover of $\inf$ is groupable
as a cover of $Y$.

Let $\cU = \seq{O_n}$ where for each $n$,
$$O_n = \{a\in\inf : n\in a\}.$$
Then $\cU$ is a large open cover of $\inf$.
Thus, there exists a partition $\seq{\cF_n}$ of $\cU$ into finite sets
such that for each $y\in Y$, $y\in\cup\cF_n$ for
all but finitely many $n$.
For each $n$ set $F_n = \{m : O_m\in\cF_n\}$.
Then $\seq{F_n}$ is a partition of $\N$ into finite sets.
For each $y\in Y$ and for all but finitely many $n$,
there exists $k$ such that $y\in O_k\in\cF_n$, that is,
$k\in y\cap F_n$, therefore $y\cap F_n\neq\emptyset$.
\end{proof}

\begin{rem}
A strengthening of the Hurewicz property for $X$,
considering \emph{countable Borel} covers instead of
open covers, was given the following simple characterization in \cite{CBC}:
\begin{quote}
For each sequence $\seq{\cU_n}$ of countable (large) Borel covers of $X$,
there exist elements $U_n\in\cU_n$, $n\in\N$, such that $X\sbst\Union_n\bigcap_{m>n}U_m$.
\end{quote}
(Note that the analogous equivalence for the open case does not
hold \cite{coc2}.) Using the same proof as in
Theorem~\ref{HureChar}, we get that this property is also
equivalent to requiring that every countable large Borel cover of
$X$ is groupable.
\end{rem}

Forgetting about the topology and considering only countable covers,
we get the following characterization of the minimal cardinality $\b$ of
an unbounded family in the Baire space $\NN$.
For a cardinal $\kappa$, denote by $\Lambda_\kappa$ (respectively, $\Lambda^{\gr}_\kappa$)
the collection of countable large (respectively, groupable) covers of $\kappa$.
\begin{cor}\label{bchar}
For an infinite cardinal $\kappa$, the following are equivalent:
\be
\i $\kappa<\b$,
\i Each subset of $\NNup$ of cardinality $\kappa$ admits a slalom,
\i For each family $Y\sbst\inf$ of cardinality $\kappa$,
there exists a partition $\seq{F_n}$ of $\N$ into
finite sets such that for each $y\in Y$,
$y\cap F_n\neq\emptyset$ for all but finitely many $n$; and
\i $\binom{\Lambda_\kappa}{\Lambda^{\gr}_\kappa}$ holds
(i.e., every countable large cover of $\kappa$ is groupable).
\ee
\end{cor}

\begin{rem}
The underlying combinatorics in this paper is similar to that appearing
in Bartoszy\'nski's characterization of $\add(\Null)$ (the minimal cardinality
of a family of measure zero sets whose union is not a measure zero set)
\cite{tomek}.
The equivalence $(1\Iff 2)$ in Corollary~\ref{bchar} is folklore.
The equivalence of these with $(4)$ seems to be new.

The only other covering property we know of which enjoys
the possibility of considering subcovers of a given cover instead of selecting
from a given sequence of covers is the Gerlits-Nagy $\gamma$-property.
The proof for this fact is much easier -- see \cite{GN}.
\end{rem}

\end{document}